# Combinatory completeness in partial groupoids

Piet Rodenburg
Instituut voor Informatica
Universiteit van Amsterdam,
december 2022

*Abstract*. I characterize the combinatorially complete pargoids (partial applicative systems) by expandability with two constants that satisfy the well-known identities. An example shows that this class contains more than just the reducts of partial combinatory algebras.

*Note*. A Dutch version of this document exists, entitled *Combinatorische volledigheid in partiële groepoïden*.

## §1. Introduction

The superficial reader of a recent paper by Shafer and Terwijn [ST] might get the impression that Feferman proved in [F] that in partial groupoids combinatory completeness — a notion the authors only loosely circumscribe — is equivalent with the presence of two elements with certain properties. However, who looks for the proof in [F], will not find it. Combinatory completeness evidently implies the existence of elements $s$ and $k$ such that $sxyz \simeq xz(yz)$ and $kxy = x$; but must $sxy$ exist for all $x$ and $y$? Feferman only states, in §3.3(1), that it is reasonable on philosophical grounds that an abstraction term $\lambda^* x.\mathbf{t}$ always has a value, and that the stipulation $sxy\!\downarrow$, in the context of the other combinator axioms, suffices for that. As it is, his minimalist treatment of possibly undefined terms — as Renardel and Troelstra pointed out already in their review [RT] — is not a good basis for a mathematical proof.

## §2. Preliminaries

I denote the transitive closure of a binary relation $R$ by $R^+$. In particular $\to^+$ is the transitive closure of the rewrite relation $\to$ (§4).

I shall follow Descartes' principle (*je pense, donc je suis*, [D]) that an atomic assertion implies the existence of its subject. In particular $M = M$ is a way of asserting that $M$ exists. More concisely I shall write $M\!\downarrow$ for '$M$ exists' and $M\!\uparrow$ for '$M$ does not exist'. The definition of the Kleene-equality $\simeq$ is:

$$M \simeq N \;:\Leftrightarrow\; M\!\downarrow \vee N\!\downarrow \to M = N.$$

A derivative of the *Cogito-Principle* stated above is the *strictness principle*: the components of a composite expression that refers to an existing object refer to existing objects. Indeed we may understand an identity $M(N) = P$ — as long as $N$ actually occurs in it — as an atomic assertion with subject $N$. In an arithmetical example: if

$$x \div (y \div z) = 6,$$

then $y \div z$ exists. (And hence $z$ is not zero.)

For equality symbol in formulas I use $\approx$. I abbreviate the formula $\mathbf{t} \approx \mathbf{t}$ to $\mathbf{t}$; the formal definition of Kleene-equality then becomes

$$\mathbf{s} \simeq \mathbf{t} \;:=\; \mathbf{s} \vee \mathbf{t} \to \mathbf{s} \approx \mathbf{t}.$$

I use Hoare's notation for simple case distinctions:

$$M \triangleleft A \triangleright N$$

is $M$ if $A$ is true, and $N$ if $A$ is false. It is not required that all cases are provided for: $M \triangleleft A$ is $M$ if $A$ is true, and undefined otherwise; and $A \triangleright N \simeq N \triangleleft \neg A$.

For algebraic notation I follow [ALV]. In particular I take a partial $n$-ary operation over a set $A$ to be a subset of $A \times A^n$. As usual a nullary operation $\{\langle a, \emptyset \rangle\}$ — a constant — will be identified with its value $a$. Since in the definitions of term operations and polynomial operations on the basis of partial fundamental operations some additional care is required, I present them briefly below.

Let $\mathbf{A}$ be a structure: a set $A$, the *domain* of $\mathbf{A}$, with a number of finitary, not necessarily total, fundamental operations over $A$. Then $\mathrm{Clo}\,\mathbf{A}$ is the clone of term operations of $\mathbf{A}$; that is, the least set of operations over $A$ that contains the fundamental operations and the projections $e_i^n \colon \langle a_0, \ldots, a_{n-1}\rangle \mapsto a_i$ of $A^n$ onto $A$, for all $i < n < \omega$, and is closed under the composition operations

$$C_m^n \colon \langle f, g_0, \ldots, g_{n-1}\rangle \mapsto f(g_0, \ldots, g_{n-1}) \qquad n, m \in \omega$$

that are defined, for $n$-ary $f$ and $m$-ary $g_0, \ldots, g_{n-1}$, by

$$f(g_0, \ldots, g_{n-1})(x_0, \ldots, x_{m-1}) \simeq f(g_0(x_0, \ldots, x_{m-1}), \ldots, g_{n-1}(x_0, \ldots, x_{m-1})).$$

For $n \in \omega$, $\mathrm{Clo}_n\,\mathbf{A}$ is the collection of $n$-ary elements of $\mathrm{Clo}\,\mathbf{A}$.

The clone $\mathrm{Pol}\,\mathbf{A}$ of polynomial operations of $\mathbf{A}$ is defined analogously; in addition to the fundamental operations it contains all constants: for each $a \in A$ the nullary operation with value $a$. For $n \in \omega$, $\mathrm{Pol}_n\,\mathbf{A}$ is the collection of $n$-ary elements of $\mathrm{Pol}\,\mathbf{A}$. In particular $\mathrm{Pol}_n\,\mathbf{A}$ contains the $n$-ary constant functions — the $n$-ary function with constant value $a$ is $C_n^0(a)$.

Let $\mathcal{L}$ be an algebraic type: a set of operation symbols of given finite arity. I shall denote the set of terms of type $\mathcal{L}$ in variables from a given set $X$ by $T_\mathcal{L}(X)$. In a structure $\mathbf{A}$ of type $\mathcal{L}$ — that specifies for every symbol $Q \in \mathcal{L}$ an operation $Q^\mathbf{A}$ over the domain $A$ of the arity proper to $Q$ — we interpret the elements of $T_\mathcal{L}(X)$ as partial functions from the set $A^X$ of assignments to $X$ in $A$, as follows:

if $\mathbf{t} = x \in X$, then $\mathbf{t}^\mathbf{A}$ is the projection $e_x^X \colon a \mapsto a(x)$;

if $\mathbf{t} = Q\mathbf{s}_1 \ldots \mathbf{s}_n$, then $\mathbf{t}^\mathbf{A} = Q^\mathbf{A}(\mathbf{s}_1^\mathbf{A}, \ldots, \mathbf{s}_n^\mathbf{A})$.

If $X = n\ (= \{0, \ldots, n-1\})$, the operations $\mathbf{t}^\mathbf{A}$ are the elements of $\mathrm{Clo}_n\,\mathbf{A}$.

Let $\mathbf{A}_A$ be the expansion of $\mathbf{A}$ in which every element of the domain $A$ is denoted by its own special constant symbol. We denote the type of $\mathbf{A}_A$ by $\mathcal{L} \cup A$. Since $\mathrm{Pol}\,\mathbf{A} = \mathrm{Clo}\,\mathbf{A}_A$, the elements of $\mathrm{Pol}_n\,\mathbf{A}$ are the operations $\mathbf{p}^\mathbf{A}$ for $\mathbf{p} \in T_{\mathcal{L} \cup A}(n)$.

Because a term in variables from $X$ also is a term in variables from any superset $Y$, the notation $\mathbf{t}^\mathbf{A}$ may be ambiguous. If necessary I shall write $\mathbf{t}_X^\mathbf{A}$ for the partial function of assignments to $X$. The least $X$ for which $\mathbf{t}_X^\mathbf{A}$ may be nonempty is $\mathrm{Var}\,\mathbf{t}$, the set of variables that occur in $\mathbf{t}$.

Usually in concrete cases, such as that of the product in a pargoid, we shall trust on the context for the distinction between an operation symbol $Q$ and an operation $Q^\mathbf{A}$, and suppress superscripts.

Let $\mathbf{A}$ and $\mathbf{B}$ be structures of the same type $\mathcal{L}$. Then $\mathbf{A}$ is a *relative substructure* of $\mathbf{B}$ if $A \subseteq B$ and for each $n$ we have, for every $n$-ary operation symbol $Q \in \mathcal{L}$ and all $a_1, \ldots, a_n \in A$:

$$Q^{\mathbf{A}}(a_1,\ldots,a_n) \simeq Q^{\mathbf{B}}(a_1,\ldots,a_n) \triangleleft Q^{\mathbf{B}}(a_1,\ldots,a_n) \in A.$$

## §3. Pargoids and combinatory completeness

A *partial applicative structure* [B], or *partial groupoid*, abbreviated *pargoid* [LE], is a nonempty set $A$ with a binairy operation · (application, product) that need not be definied for all pairs in $A \times A$. A pargoid is *total*, or a *groupoid* (or *applicative structure*) if its product operation is total.

An element $a$ of a pargoid $\mathbf{A} = \langle A, \cdot \rangle$ is *left passive* if for no $x \in A$ the product $a \cdot x$ is defined.

The operation symbol · is often suppressed. There are three kinds of contexts in which we *do* write the operation dot: when the operation has a certain emphasis, as in definitions; to highlight one of the applications, for example to save on parentheses, when we write $x \cdot yz$ instead of $x(yz)$; and when some factors are represented by complex expressions. We let the product associate to the left: $xyz = (xy)z$.

A non-total product operation induces totally undefined polynomial operations:

PROPOSITION 1. Let $\mathbf{A} = \langle A, \cdot \rangle$ be a non-total pargoid. Then $\emptyset \in \text{Pol}_n \mathbf{A}$ for all $n \geq 0$.

PROOF. Suppose $a_1 a_2 \uparrow$. Then $C_0^2(\cdot, a_1, a_2)$ is the empty nullary operation; and $C_n^0(C_0^2(\cdot, a_1, a_2))$ is $n$-ary and empty. ⊠

DEFINITION 2. A pargoid $\mathbf{A}$ is *combinatorially complete* if for all $n \geq 1$, for every $n$-ary polynomial operation $p$ of $\mathbf{A}$, an element $a \in A$ exists such that
for all $x_1,\ldots,x_n \in A$, $\quad ax_1\ldots x_n \simeq p(x_1,\ldots,x_n)$.

REMARK. At $n = 0$ a problem appears. If $p \neq \emptyset$, then $p$ itself is a suitable $a$. But if $\mathbf{A}$ is not total, then there also is a nullary polynomial operation without a value; and this operation is of course *not* equivalent to an element of $A$.

Definition 2 generalizes the usual definition for groupoids:

PROPOSITION 3. A groupoid $\mathbf{A}$ is combinatorially complete if and only if for all $n \geq 0$, for every $n$-ary polynomial operation $p$ of $\mathbf{A}$, some $a \in A$ exists such that for all $x_1,\ldots,x_n \in A$, $ax_1\ldots x_n = p(x_1,\ldots,x_n)$.

PROOF. Observe that in a groupoid all polynomial operations are total; hence for polynomial operations $p$ and $q$ and $x_1,\ldots,x_n \in A$:

(∗) $\qquad p(x_1,\ldots,x_n) \simeq q(x_1,\ldots,x_n) \Leftrightarrow p(x_1,\ldots,x_n) = q(x_1,\ldots,x_n)$.

($\Rightarrow$) Let $p$ be a nonempty nullary polynomial operation. Then $p$ is a singleton $\{\langle a, \emptyset \rangle\}$, for some $a \in A$; so $a = p(\emptyset)$.
($\Leftarrow$) Immediately from (∗). ⊠

## §4. Partial combinatory algebras

A *semicombinatory algebra*, abbreviated *sca*, is a pargoid $\mathbf{A}$ in which two constant symbols $s$ and $k$ are interpreted so that

$$\mathbf{A} \models (1)\ sxyz \simeq xz(yz)\ \&\ (2)\ kxy \approx x.$$

An sca is a *partial combinatory algebra*, abbreviated *pca*, if moreover

$$\mathbf{A} \models (0)\ sxy.$$

A pca is *total*, or a *combinatory algebra*, if its product operation is total.

LEMMA 1 (Bethke). Every non-total pca contains a left passive element.

PROOF. Let $\mathbf{A} = \langle A, \cdot, s, k \rangle$ be a pca, and $a$, $b$ elements of $A$ whose product is undefined. Since $kaa = a$ and $kbb = b$, the Cogito-Principle implies that $ka$ and $kb$ exist; hence by (0), $s(ka)(kb)$ exists. By (1) and (2) we now have for every $x \in A$:

$$s(ka)(kb)x \simeq kax(kbx) \simeq ab,$$

so $s(ka)(kb)x\uparrow$. ⊠

THEOREM 2. A pargoid is combinatorially complete if and only if it can be expanded to a semicombinatory algebra, and either is total, or contains a left passive element.

PROOF. ($\Rightarrow$) Let $\mathbf{A}$ be a combinatorially complete pargoid. Combinatory completeness implies that $\mathbf{A}$ has elements $s$ and $k$ that satisfy

$$(1)\ sxyz \simeq xz(yz)\ \&\ (2a)\ kxy \simeq x.$$

Since $x$ varies over existing elements of $\mathbf{A}$, (2a) is equivalent to

$$(2)\ kxy \approx x.$$

If $\mathbf{A}$ is not total, we have $\emptyset \in \text{Pol}_1 \mathbf{A}$ by §3, Proposition 1. By combinatory completeness there must then be $a \in A$ such that for all $x \in A$

$$ax \simeq \emptyset(x),$$

so $ax\uparrow$.

($\Leftarrow$) Let $\mathbf{A} = \langle A, \cdot, s, k \rangle$ be a semicombinatory algebra. Let $C$ be the type $\{\cdot, s, k\}$, and $T_{C \cup A}(Y)$ the set of formal $\mathbf{A}$-polynomials in a countably infinite set $Y$ of variables. We define in the usual way for each $x \in Y$ an operation $\lambda^* x$ that maps $T_{C \cup A}(Y)$ to $T_{C \cup A}(Y \setminus \{x\})$: $\lambda^* x.\mathbf{p}$ is

$\quad \mathbf{i} := skk$ if $\mathbf{p} = x$,

$\quad k\mathbf{p}$ if $x$ does not occur in $\mathbf{p}$,

$\quad s(\lambda^* x.\mathbf{r})(\lambda^* x.\mathbf{q})$ otherwise, if $\mathbf{p} = \mathbf{rq}$.

A simple induction over $\mathbf{p}$ shows that $x$ does not occur in $\lambda^* x.\mathbf{p}$. Next we prove, again by induction over $\mathbf{p}$, that for all $b \in A$ en $c: Y \setminus \{x\} \to A$,

$$(\lambda^* x.\mathbf{p})^{\mathbf{A}}_{Y \setminus \{x\}}(c) \cdot b \simeq \mathbf{p}^{\mathbf{A}}_Y(c \cup \{\langle b, x \rangle\}) \qquad (\lambda^*).$$

We simplify the expression on the righthand side to $\mathbf{p}^{\mathbf{A}}_Y(c, b)$.

1° $\mathbf{p} = x$. It is to be proved that $skkb = b$, for any $b \in A$. By (1),

$$skkb \simeq kb(kb);$$

by (2), $kbb$ exists, and hence so does $kb$; by (2) once more, $kb(kb) = b$.

2° $x$ does not occur in $\mathbf{p}$. Take $c \in A^{Y \setminus \{x\}}$. We must show that

$$k \cdot \mathbf{p}^{\mathbf{A}}_{Y \setminus \{x\}}(c) \cdot b \simeq \mathbf{p}^{\mathbf{A}}_Y(c, b), \qquad \text{for each } b \in A.$$

Observe that $\mathbf{p}^{\mathbf{A}}_{Y\setminus\{x\}}(c) \simeq \mathbf{p}^{\mathbf{A}}_Y(c, b)$. If $\mathbf{p}^{\mathbf{A}}_{Y\setminus\{x\}}(c)$ exists, then
$$k \cdot \mathbf{p}^{\mathbf{A}}_{Y\setminus\{x\}}(c) \cdot b = \mathbf{p}^{\mathbf{A}}_{Y\setminus\{x\}}(c)$$
by (2). If on the other hand $\mathbf{p}^{\mathbf{A}}_{Y\setminus\{x\}}(c)$ does not exist, then $k \cdot \mathbf{p}^{\mathbf{A}}_{Y\setminus\{x\}}(c) \cdot b$ does not exist either, by the strictness principle.

3° $\mathbf{p} = \mathbf{rq}$. Take $c \in A^{Y\setminus\{x\}}$. We are to prove that
$$s \cdot (\lambda^* x.\mathbf{r})^{\mathbf{A}}_{Y\setminus\{x\}}(c) \cdot (\lambda^* x.\mathbf{q})^{\mathbf{A}}_{Y\setminus\{x\}}(c) \cdot b \simeq \mathbf{r}^{\mathbf{A}}_Y(c, b) \cdot \mathbf{q}^{\mathbf{A}}_Y(c, b), \quad \text{for each } b \in A.$$
If $\mathbf{r}^{\mathbf{A}}_Y(c, b) \cdot \mathbf{q}^{\mathbf{A}}_Y(c, b)$ exists, then so do $\mathbf{r}^{\mathbf{A}}_Y(c, b)$ and $\mathbf{q}^{\mathbf{A}}_Y(c, b)$, by the strictness principle; by induction hypothesis $(\lambda^* x.\mathbf{r})^{\mathbf{A}}_{Y\setminus\{x\}}(c)$ and $(\lambda^* x.\mathbf{q})^{\mathbf{A}}_{Y\setminus\{x\}}(c)$ then exist as well. Hence we may apply (1). For the converse direction, suppose
$$s \cdot (\lambda^* x.\mathbf{r})^{\mathbf{A}}_{Y\setminus\{x\}}(c) \cdot (\lambda^* x.\mathbf{q})^{\mathbf{A}}_{Y\setminus\{x\}}(c) \cdot b \text{ exists.}$$

Then
$$\begin{aligned}
s \cdot (\lambda^* x.\mathbf{r})^{\mathbf{A}}_{Y\setminus\{x\}}(c) \cdot (\lambda^* x.\mathbf{q})^{\mathbf{A}}_{Y\setminus\{x\}}(c) \cdot b &= \\
(\lambda^* x.\mathbf{r})^{\mathbf{A}}_{Y\setminus\{x\}}(c) \cdot b \cdot ((\lambda^* x.\mathbf{q})^{\mathbf{A}}_{Y\setminus\{x\}}(c) \cdot b) &\quad \text{by (1)} \\
= \mathbf{r}^{\mathbf{A}}_Y(c, b) \cdot \mathbf{q}^{\mathbf{A}}_Y(c, b) &\quad \text{by induction hypothesis} \\
= \mathbf{p}^{\mathbf{A}}_Y(c, b).
\end{aligned}$$

By the strictness principle, $(\lambda^*)$ implies:
$$\operatorname{Dom}(\lambda^* x.\mathbf{p})^{\mathbf{A}}_{Y\setminus\{x\}} \supseteq \{c \in A^{Y\setminus\{x\}} \mid \exists b \in A.\ c \cup \{\langle b, x\rangle\} \in \operatorname{Dom}\mathbf{p}^{\mathbf{A}}_Y\} \quad (\mathrm{D}\lambda^*).$$

Let $d$ be a left passive element of $\mathbf{A}$. Take $p \in \operatorname{Pol}_n \mathbf{A}$. Assume $p = \mathbf{p}^{\mathbf{A}}$, with $\mathbf{p} \in T_{C \cup A}(x_1, \ldots, x_n)$; define
$$a := d \triangleleft \operatorname{Dom} p = \emptyset \triangleright (\lambda^* x_1. \ldots \lambda^* x_n.\mathbf{p})^{\mathbf{A}}_\emptyset.$$
If $\operatorname{Dom} p = \emptyset$, then $ab_1 \ldots b_n$ and $p(b_1, \ldots, b_n)$ equally diverge for all $b_1, \ldots, b_n \in A$. Now suppose $\operatorname{Dom} p \neq \emptyset$. Then we conclude, in $n$ steps, using $(\mathrm{D}\lambda^*)$, that $(\lambda^* x_1. \ldots \lambda^* x_n.\mathbf{p})^{\mathbf{A}}_\emptyset\downarrow$. Next we calculate
$$\begin{aligned}
ab_1 \ldots b_n &\simeq (\lambda^* x_1. \ldots \lambda^* x_n.\mathbf{p})^{\mathbf{A}} b_1 \ldots b_n &&\text{by definition} \\
&\simeq (\lambda^* x_2. \ldots \lambda^* x_n.\mathbf{p}(b_1, x_2, \ldots, x_n))^{\mathbf{A}} b_2 \ldots b_n &&\text{by } (\lambda^*) \\
&\simeq (\lambda^* x_3. \ldots \lambda^* x_n.\mathbf{p}(b_1, b_2, x_3, \ldots, x_n))^{\mathbf{A}} b_3 \ldots b_n &&\text{likewise} \\
&\simeq \ldots\ldots \\
&\simeq \mathbf{p}(b_1, \ldots, b_n)^{\mathbf{A}} = p(b_1, \ldots, b_n).
\end{aligned}$$
We conclude that $\mathbf{A}$ is combinatorially complete.

If $\mathbf{A}$ does not have a left passive element, then by Lemma 1, $\mathbf{A}$ is total. Then $\mathbf{A}$ has no empty polynomial operations, and it suffices to define
$$a := (\lambda^* x_1. \ldots \lambda^* x_n.\mathbf{p})^{\mathbf{A}}_\emptyset$$
to show that $\mathbf{A}$ is combinatorially complete. ⊠

CORLLARY 3. Partial combinatory algebras are combinatorially complete.

PROOF. Partial combinatory algebras are semicombinatory algebras, and the non-total ones contain a left passive element, by Lemma 1. ⊠

**Example 4**

Let $C = \{\cdot, s, k\}$ be the type of combinatory algebras, $X$ an arbitrary set of variables, and $N$ the set of normal forms in $T_C(X)$ of the term rewriting system CL :=

$$sxyz \to xz(yz),$$
$$kxy \to x.$$

Define the product operation $*$ over $N$ by

$$\mathbf{m} * \mathbf{n} := \text{the normal form of } \mathbf{mn}.$$

(The system CL is confluent; so normal forms are unique. Cf. [Te].)

The structure $\mathbf{N} = \langle N, *, s, k \rangle$ is a pca: the rewrite rules of CL correspond to the laws (1) and (2), and if $\mathbf{m}$ and $\mathbf{n}$ are normal forms of CL, then $s\mathbf{mn}$ is a normal form as well. Hence Corollary 3 implies that $\mathbf{N}$ (or rather its pargoid reduct) is combinatorially complete.

The converse of Corollary 3 is

(†)  Every combinatorially complete pargoid is expandable to a pca.

By a variation on Example 4 we will show that (†) is *not* true.

**Example 5**

Let $\mathbf{N}$ be as in the previous example; only assume that $X$ is infinite. Let $\omega := s\mathbf{ii}$. The term $\omega$ belongs to $N$, but $\omega * \omega$ does not exist. Define: $\mathbf{d} := s(k\omega)(k\omega)$. Then $\mathbf{d}$ is a left passive element of $\mathbf{N}$ (by the proof of Lemma 1). Now let

$$L := \{\mathbf{n} \in N \mid (\mathbf{n} * x)\!\downarrow \text{ for } x \in X \setminus \text{Var}\mathbf{n}\},$$

and $N' = L \cup \{\mathbf{d}\}$. Observe that $L$ contains all terms in which the constant symbols $s$ and $k$ do not occur. Let $\mathbf{N}'$ be the relative substructure of $\mathbf{N}$ with domain $N'$.

PROPOSITION 5.1. $\mathbf{N}'$ is an sca.

PROOF. Let $\mathbf{p}$ and $\mathbf{q}$ be elements of $T_C(X)$ such that $\mathbf{p} \to \mathbf{q}$ in CL. Then either $\mathbf{p}$ and $\mathbf{q}$ have the same normal form, and this normal form belongs to $N'$, or neither has a normal form in $N'$. So $\mathbf{N}' \models \mathbf{p} \simeq \mathbf{q}$; and if $\mathbf{q} \in X, \mathbf{N}' \models \mathbf{p} \approx \mathbf{q}$. Thus $\mathbf{N}'$ satisfies the equalities (1) and (2).  □

COROLLARY 5.2. $\mathbf{N}'$ is combinatorially complete.

PROOF. Combine Theorem 2 with the above Proposition and the fact that a left passive element of a pargoid $\mathbf{P}$ is also left passive in relative substructures of $\mathbf{P}$.  □

The CL-normal form $\omega$ belongs to $L$, for $\omega * x = xx$. Hence also $k\omega \in L$. Likewise $s\omega\mathbf{i} \in L$: $s\omega\mathbf{i} * x = xxx$; and $k(s\omega\mathbf{i}) \in L$. Now observe that $s\omega\mathbf{i} * \omega\!\uparrow$:

$$s\omega\mathbf{i}\omega \to \omega\omega(\mathbf{i}\omega) \to^+ \omega\omega\omega \to^+ \omega\omega\omega.$$

Therefore the normal form $s(k(s\omega\mathbf{i}))(k\omega)$ is left passive in $\mathbf{N}$ (cf. the proof of Lemma 1); but $s(k(s\omega\mathbf{i}))(k\omega) \neq \mathbf{d}$, so $s(k(s\omega\mathbf{i}))(k\omega) \notin N'$. Consequently $\mathbf{N}' \not\models (0)$; so $\mathbf{N}'$ is not a pca.

This does not yet refute (†), however, since formally combinatory completeness is a property of the pargoid $\mathbf{N}'\!\upharpoonright\{\cdot\}$.

PROPOSITION 5.3. There are no $\mathbf{t}, \mathbf{u} \in N'$ such that $\langle N', *, \mathbf{t}^\mathbf{N}, \mathbf{u}^\mathbf{N}\rangle$ is a pca.

Proof. We shall prove that there is no $\mathbf{t} \in N'$ such that

(3) $\qquad\qquad \mathbf{N}' \models \mathbf{t}xy \;\&\; \mathbf{t}xyz \simeq xz(yz).$

Suppose $\mathbf{t}$ satisfies (3). Then $\mathbf{N}' \models \mathbf{t}(k(s\omega\mathbf{i}))(k\omega)$. Take three variables $x, y$ and $z$ that do not belong to Var$\mathbf{t}$. By (3) we have, since $xz(yz)\downarrow$,

$$\mathbf{t}xyz \to^+ xz(yz);$$

substitution of $k(s\omega\mathbf{i})$ for $x$ and $k\omega$ for $y$ in the reduction results in

$$\mathbf{t}(k(s\omega\mathbf{i}))(k\omega)z \to^+ k(s\omega\mathbf{i})z(k\omega z).$$

But $k(s\omega\mathbf{i})z(k\omega z)$ reduces to $s\omega\mathbf{i}\omega$, hence does not have a normal form. So

$$(\mathbf{t}(k(s\omega\mathbf{i}))(k\omega))_X^{\mathbf{N}'}(1_X) \notin L;$$

hence the normal form of $\mathbf{t}(k(s\omega\mathbf{i}))(k\omega)$ must be $\mathbf{d}$. We conclude that on the one hand $\mathbf{t}(k(s\omega\mathbf{i}))(k\omega)z$ reduces to $s\omega\mathbf{i}\omega$, and on the other to $\mathbf{d}z$, and a fortiori to $\omega\omega$. Since CL is confluent, $\omega\omega\omega$ and $\omega\omega$ should then have a common reduct. But they clearly have not. ☒

**Example 6**

Turing described in [Tu] a notion of computing machine, and argued that it covers every conceivable algorithm. There exist bijective enumerations of the Turing-machines and their computations such that the predicate $Texy$,

"on input $x$ Turing-machine $e$ terminates after computation $y$",

is primitive recursive, and the function $U$ that reads off the result, if there is one, from the number representing a computation, computable. The computable function with index $e$ is then defined by

$$\varphi_e(x) \simeq U(\mu y.Texy).$$

By the recursion theorem there are total computable functions $f$ and $g$ such that for all $x, y, z \in \mathbb{N}$,

$$\varphi_{f(x)}(y) \simeq x \quad \text{and} \quad \varphi_{g(x,y)}(z) \simeq \varphi_{\varphi_x(z)}(\varphi_y(z));$$

and a total computable function $h$ such that for all $x$ and $y \in \mathbb{N}$, $\varphi_{h(x)}(y) = g(x, y)$.

Now let $s$ be an index of $h$, and $k$ an index of $f$. Define a product operation over $\mathbb{N}$ by

$$x \cdot y \simeq \varphi_x(y).$$

Then

$$k \cdot x = f(x) \quad \text{en} \quad s \cdot x \cdot y = g(x, y);$$

a fortiori the laws (1) and (2) holds. So $\langle \mathbb{N}, \cdot, s, k \rangle$ is a partial combinatory algebra.

## §5. Conclusion

The foregoing may be viewed as yet another illustration of the complexity of reasoning with partial operations.

The definition of polynomial operations is in my opinion as general as possible. A less general definition, by the way, would increase the class of combinatorially complete pargoids, so that there would not be fewer counterexamples to statement (†) above Example 5.

Feferman's argument applies to the computability theoretic pargoids studied by Shafer and Terwijn in [ST] — Example 6 describes the simplest of them.

Thus the function *g* in Example 6 transforms two numbers into an algorithmic code, and if we do not expect too much of the codes, such a transformation is always possible. There is, however, no intrinsic reason why a combinatorially complete pargoid should allow a division between codes on the one hand and risky computations on the other.

*Acknowledgement*. A discussion with Henk Barendregt about an early version of this paper has led to a considerable improvement of the presentation.